\documentclass[12pt]{amsart}
\usepackage[utf8]{inputenc}

\usepackage[margin=1.2in]{geometry}                
\usepackage{amsmath,amsthm,amssymb}
\usepackage{graphicx}
\usepackage{comment}
\usepackage{epstopdf}
\usepackage{tikz}
\usepackage{tikz-cd} 
\usepackage{xcolor}
\usepackage[colorlinks=true, pdfstartview=FitV, linkcolor=blue, citecolor=blue, urlcolor=blue]{hyperref}
\usepackage[nameinlink]{cleveref}
    \crefname{conj}{conjecture}{conjectures}
    \crefname{algocfline}{algorithm}{algorithms}
    
\usepackage{pgfplots, multicol, float, mathtools, extarrows}
\usepgfplotslibrary{fillbetween}
\usepackage[linesnumbered,ruled,vlined]{algorithm2e}

\makeatletter  
\@namedef{subjclassname@2020}{%
  \textup{2020} Mathematics Subject Classification}
\makeatother

\newcommand{\N}{{\mathbb N}}

\newcommand{\R}{{\mathbb R}}
\newcommand{\Q}{{\mathbb Q}}
\newcommand{\Z}{{\mathbb Z}}

\newcommand{\ba}{{\mathbf a}}

\newcommand{\bc}{{\mathbf c}}

\newcommand{\bp}{{\mathbf p}}

\newcommand{\bx}{{\mathbf x}}
\newcommand{\by}{{\mathbf y}}
\newcommand{\bz}{{\mathbf z}}

\def\Ass{\operatorname{Ass}}

\def\conv{\operatorname{convex\ hull}}

\def\het{\operatorname{ht}}

\def\ov{\overline}

\newtheorem{thm}{Theorem}[section]

\newtheorem*{introthm*}{Theorem}
\newtheorem{conj}[thm]{Conjecture}
\newtheorem{cor}[thm]{Corollary}
\newtheorem{lem}[thm]{Lemma}
\newtheorem{prop}[thm]{Proposition}

\newtheorem{quest}[thm]{Question}

\theoremstyle{definition}
\newtheorem{defn}[thm]{Definition}

\newtheorem{ex}[thm]{Example}

\theoremstyle{remark}
\newtheorem{rem}[thm]{Remark}
\newtheorem{notation}[thm]{Notation}

\numberwithin{equation}{section}  

\title{Consequences of the packing problem}
\author[Polymath 2020, Monomials, Convex Bodies, and Optimization Team]{Hrishikesh Bodas, Benjamin Drabkin, Caleb Fong, Su Jin, Justin Kim, Wenxuan Li, Alexandra Seceleanu, Tingting Tang, Brendan Williams}
\thanks{The second author was supported by the NSF RTG grant in algebra and combinatorics at the University of Minnesota  DMS--1745638. The seventh author was supported by NSF DMS--1601024. This work was completed in the framework of the 2020 Polymath program \url{https://geometrynyc.wixsite.com/polymathreu}.}
\address{Carnegie Mellon University}
\email{hbodas@andrew.cmu.edu}
\address{University of Nebraska--Lincoln}
\email{benjamin.drabkin@huskers.unl.edu}
\address{University of St Andrews}
\email{cjxf@st-andrews.ac.uk}
\address{University of Illinois, Urbana-Champaign}
\email{sujin2@illinois.edu}
\address{Vanderbilt University}
\email{hanbin919@gmail.com}
\address{University of California, Santa Barbara}
\email{wenxuanli@ucsb.edu}
\address{University of Nebraska--Lincoln}
\email{aseceleanu@unl.edu}
\address{San Diego State University}
\email{ttang2@sdsu.edu}
\address{University of Michigan--Dearborn}
\email{brendwil@umich.edu}

\begin{document}

\begin{abstract}
We study several consequences of the packing problem, a conjecture from combinatorial optimization, using algebraic invariants of square-free monomial ideals. While the packing problem is currently unresolved, we successfully settle the validity of its consequences. Our work prompts additional questions and conjectures, which are presented together with their motivation.
\end{abstract}

\maketitle


\section{Introduction}

The packing problem  introduced by Conforti and Cornu\'ejols \cite{conforti1990decomposition} is a conjecture originating from combinatorial optimization in the context of max-flow min-cut properties. It has been brought into commutative algebra through the inspiring paper \cite{GVV} of Gitler, Valencia and Villarreal. A comprehensive account of this problem from an algebraic perspective can be found in the monograph \cite[\S 14.3]{MonAlgBook} as well as in the surveys \cite{francisco2013powers} and \cite{dao2017symbolic}. We state an algebraic version of this problem  in \Cref{conj:packingalg} as well as a more combinatorial version that is closer to the roots of the problem in combinatorial optimization in \Cref{conj:packing}. 

Our contribution in this paper is to study several consequences of the packing problem. Since the packing problem is at the moment of this writing still a conjecture, we believe that it is useful to approach it gradually by establishing the truth for at least some of its consequences. We also make connections between the circle of ideas related to the packing problem and Alexander duality for square-free monomial ideals.

To state the packing problem and our main results, we need to introduce some of the main characters of this writing.  A square-free monomial ideal is an ideal generated by square-free monomials in a polynomial ring. This class of ideals encodes numerous combinatorial objects, chief among which are simplicial complexes and hypergraphs. In this paper we focus on the correspondence between  square-free monomial ideals and hypergraphs.

A square-free monomial ideal $I$ has a unique irredundant decomposition into prime ideals which takes the form $I = P_1 \cap \cdots \cap P_s$ with $P_i = ( x_{j_1},\ldots,x_{j_{s_j}})$ for $j=1,\ldots,s$. Furthermore, one defines the height of $I$, $\het(I)$ to be the minimum of the number of minimal generators of the prime components $P_i$.
Based on the above decomposition, for each positive integer $m$ one defines the $m$-th {\em symbolic power} of $I$ to be the monomial ideal
\[I^{(m)} = P_1^m \cap \cdots \cap P_s^m .\]
Symbolic power ideals are important in algebraic geometry where they encode polynomial functions vanishing to high order on a given algebraic variety. They are also relevant in combinatorics.  For example, if $I$ is the edge ideal of a (hyper)graph, the $m$-th symbolic power of $I$ encodes the $m$-covers of the (hyper)graph.

 We single out a class of square-free monomial ideals which is important to this project. 
 \begin{defn}
 \label{def:packed}
A square-free monomial ideal $I$ is {\em K\"onig} if there is a set of  pairwise coprime monomials in $I$ of cardinality $\het(I)$. 
 
The ideal $I$ has the {\em packing property} if every ideal obtained from $I$ by setting a (possibly empty) subset of the variables equal to 0 and a disjoint (possibly empty) subset of the variables equal to 1 is K\"onig.
\end{defn}

The terminology  K\"onig is best explained by the connection to  K\"onig's theorem on bipartite graphs; see the discussion preceding \Cref{thm:Konig} and the terminology packed is explained by the relationship to edge packings in hypergraphs. See \Cref{conj:packing} for a combinatorial formulation of the packing problem which clarifies this perspective. 

\begin{conj}[The packing problem -- {\cite[Conjecture 3.10]{GVV}, \cite[Theorem 4.6]{GRV}}]
\label{conj:packingalg}
 The symbolic and ordinary powers of a square-free monomial ideal $I$ coincide, i.e. $I^{(m)}=I^m$ for all positive integers $m$, if and only if $I$ has the packing property.
\end{conj}

While the direct implication of \Cref{conj:packingalg} is known to hold, cf.~ \cite[Theorem 4.6, Corollary 4.14]{GRV}; see also \cite[p.~422]{dao2017symbolic}, the converse implication is at the time of this writing a long standing conjecture. For the case of graphs, \Cref{conj:packing} holds by \cite[Proposition 4.27, Theorem 4.6]{GRV}.  Previous work on the packing problem includes \cite{CGMpacking, Corn,  HaMorey,  DR, MoreyVillarreal, MontanoNunez, AlilooeeBanerjee}. 

Our work establishes that three consequences of the converse implication in the packing problem hold. The first consequence gives a numerical shadow of the equality of the ordinary and symbolic powers of an ideal in the form of an equality between the initial degree and the Waldschmidt constant.  The initial degree $\alpha(I)$ of a homogeneous ideal $I$ is  the least degree of a nonzero element of $I$. The Waldschmidt constant of $I$ can be viewed as an asymptotic initial degree for the family of symbolic powers of $I$. This invariant is defined as $\widehat{\alpha}(I)=\lim_{m\to \infty} \frac{\alpha(I^{(m)})}{m}$; see \Cref{def:Waldschmidt} for details.

The second consequence gives a shadow
of the equality of the ordinary and symbolic powers of an ideal in convex geometric terms. In detail, there are two convex bodies that can be associated to the families of ordinary and  symbolic powers of a monomial ideal $I$ respectively, see \cite{Waldschmidtpaper}. These are termed the Newton polyhedron of $I$, $NP(I)$, and the symbolic polyhedron of $I$, $SP(I)$ cf. \Cref{def:NP} and \Cref{def:SP}. We show that these two polyhedra are equal for ideals which have the packing property. Equivalently linear programs having these two convex bodies as feasible sets have the same solutions.

Our main results on consequences of the packing problem are summarized below:

\begin{introthm*}[\Cref{thm:packing_and_SP}, \Cref{conj:packing_and_LP}, \Cref{conj:packing_and_waldschmidt}]
If $I$ is a square-free monomial ideal which satisfies the packing property then there are equalities $\alpha(I)=\widehat{\alpha}(I)$ and $NP(I)=SP(I)$, as predicted by the packing problem. Moreover the optimal solution for any linear program with feasible set $NP(I)$ coincides with the optimal solution for the linear program with the same objective function and feasible set $SP(I)$.
\end{introthm*}

We also study the relationship between the packing property and Alexander duality, with the following conclusion.

\begin{introthm*}[\Cref{prop:equidim}]
Let $I$ be an equidimensional square-free monomial ideals $I$ such that $I^\vee$ is also equidimensional.  Then $I$ and  $I^\vee$ satisfy the packing property simultaneously, that is, $I$ satisfies the packing property if and only if $I^\vee$ does.
\end{introthm*}

Our paper is organized as follows: \cref{s:hyp} provides a dictionary between square-free monomial ideals and hypergraphs, presents several combinatorial optimization invariants of hypergraphs and restates the packing problem in combinatorial language. \Cref{s:alginv} introduces several  convex bodies and combinatorial optimization invariants for monomial ideals. Three consequences of the packing problem are introduced and proven in \cref{s:consequences} as \Cref{thm:packing_and_SP}, \Cref{conj:packing_and_LP}, \Cref{conj:packing_and_waldschmidt}. In \cref{s:reverse}, we discuss the irreversibility of the consequences of the packing problem formulated in this paper, single out the class of uniform hypergraphs  as a possible candidate for which the converses of our results may apply, and prove \Cref{prop:equidim} regarding the relationship between the packing property and Alexander duality.

\section{Square-free monomial ideals and hypergraphs}
\label{s:hyp}

In this section we present the fundamental dictionary relating square-free monomial ideals to hypergraphs, also known as clutters. We supplement this dictionary by interpreting some hypergraph and ideal theoretic invariants by means of linear optimization. An excellent reference for this theory is \cite{MonAlgBook}. We do not make any claims of novelty for the contents of this section. Much of it can be found in \cite{HaTrung}.

We denote by $\N$ the set of non negative integers and by $[n]$ the set $\{1,\ldots,n\}$.

\subsection{Square-free monomial ideals as edge ideals of hypergraphs}

An ideal of the polynomial ring $R=K[x_1,\ldots, x_n]$ with coefficients in a field $K$ is a {\it monomial ideal} if it is generated
by monomials. It is a {\it square-free
monomial ideal} if it is generated by square-free monomials,
i.e., every generator has the form $x_{i_1}\cdots x_{i_t}$
with $i_j \in [n]$.
A square-free monomial ideal $I$ has a unique irredundant decomposition into prime ideals which takes the form
\[I = P_1 \cap \cdots \cap P_s ~~\mbox{with $P_i = 
( x_{j_1},\ldots,x_{j_{s_j}})$ for
$j\in[s]$}.\]
The prime ideals $P_j$ appearing in this decomposition are called the {\em associated primes} of $I$ and  form a set denoted $\Ass(I)$.
The {\em height} of a square-free monomial ideal $I$ is the least number of variables needed to generate any of its associated primes, i.e.,
\[
\het(I)=\min_{P\in \Ass(I)}{\het(P)}=\min_{j\in[s]}  s_j.
\]

A {\it hypergraph} is an ordered pair $H = (V,E)$ where
$V = \{x_1,\ldots,x_n\}$ is the set of {\it vertices}, and
$E$ consists of subsets of $V$ such that if $e_i \subseteq e_j$,
then $e_i = e_j$.  The elements of $E$ are called {\it edges}.
When the cardinality of each edge is  $|e_i|=2$, $H$ is a {\it graph}.

There is a bijective correspondence between hypergraphs $H$ on $n$ vertices and
square-free monomial ideals of $R=K[x_1,\ldots, x_n]$ by means of the following construction.

\begin{defn}
Given any hypergraph $H = (V,E)$, one can associate to
$H$ a square-free monomial ideal $I(H)$ called the {\it edge ideal}
of $H$.  Precisely, we define
\[I(H) = ( x_{i_1}x_{i_2} \cdots x_{i_t} ~|~
\{i_1,i_2,\ldots,i_t\} \in E ).\]
\end{defn}

The correspondence between hypergraphs $H$ and square-free monomial ideals $I(H)$ extends to a dictionary relating combinatorial invariants of $H$ to algebraic invariants of $I(H)$. For example, the associated primes of $I(H)$ are related to the
maximal independent sets and minimal vertex covers of the hypergraph $H$. We say that $A \subseteq V$ is
an {\it independent set} of $H$ if $e \not\subseteq A$ whenever $e \in E$.
It is {\it maximal} if it is maximal with respect to inclusion among all independent sets of $H$.

 A subset $U \subseteq V$ is a {\it vertex cover} or {\em transversal} of a hypergraph if $e \cap U \neq \varnothing$ whenever $e \in E$.
A vertex cover is {\it minimal} if it is so with respect to containment.  A {\em minimum vertex cover} is a vertex cover of smallest cardinality. Note that a minimum vertex cover is minimal, but the converse need not be true. The cardinality of any minimum vertex cover for a hypergraph $H$ is denoted $\tau(H)$ and termed the {\em transversal number} of the hypergraph $H$.

The following lemma gives a formal description of the relationship between associated primes of $I(H)$ and minimal vertex covers and maximal independent sets of $H$.
\begin{lem}
\label{vertexcovers}
Suppose that $H = (V,E)$ is a hypergraph with $E\neq \emptyset$ 
and let $I = I(H)$.
Let $I = P_1 \cap
\cdots \cap P_s$ be the irredundant prime decomposition of $I$, and
set $P'_i = \{ x_j ~|~ x_j \not\in P_i\}$ for $i\in[s]$.  Then, identifying the set of generators for each of these ideals with the set of corresponding vertices in $V$, yields
\begin{enumerate}
\item $P_1,\ldots,P_s$ are the minimal vertex covers of $H$
\item $P'_1,\ldots,P'_s$ are the maximal independent sets of $H$,
\item $\het(I(H))=\tau(H)$.
\end{enumerate}
\end{lem}

\begin{proof}
The first statement is proven in \cite[Lemma 6.3.37]{MonAlgBook}.
The last statement follows from the first and the definitions for height and $\tau(H)$.
For the second statement, any $P'_i$ is a maximal independent
set if and only if $V \setminus P_i$ is a minimal vertex cover.
We now use the first claim to finish the proof.
\end{proof}

\subsection{Linear optimization invariants of hypergraphs}
\label{s:hypinv}

The transversal number of a hypergraph introduced above can be described as the solution of an integer optimization problem. To formulate the problem, we introduce incidence matrices.
 
 \begin{defn}
 \label{def:matB}
  The {\em incidence matrix} of the hypergraph $H=(V,E)$ with $V=\{v_1,\ldots, v_n\}$ and $E=\{e_1,\ldots, e_t\}$ is the $n\times t$ matrix given by
\begin{equation}
\label{eq:matB}
B_{i,j} = \begin{cases}
1 & \mbox{if $v_i \in e_j$} \\
0 & \mbox{if $v_i \not\in e_j$.}
\end{cases}
\end{equation}
\end{defn}

The following lemma utilizes the notation $\bz=\begin{bmatrix} z_1 &\cdots &z_n\end{bmatrix}^T$ for a column vector in $\R^n$ , ${\bf 0}$ for the zero vector in $\R^n$, and  ${\bf 1}$ for the vector in $\R^n$ with all entries equal to 1. Moreover, inequalities between vectors are understood componentwise.
 
 \begin{lem}
 \label{lem:transversalopt}
The transversal number $\tau(H)$ of a hypergraph $H$  is the optimum value of the following  integer program 
\begin{equation}
\label{eq:tauLP}
\begin{tabular}{rl}
\emph{minimize} & $z_1+\cdots+z_n $\\
\emph{subject to} & $B^T{\bf z} \geq {\bf 1}$ \emph{and}  ${\bf z}\in\N^n$.
\end{tabular}
\end{equation}
\end{lem}
\begin{proof}
A vector $\bz\in \Z^n$ satisfies the inequality $B^T{\bf z} \geq {\bf 1}$ if and only if for each edge $e_i\in E$ there is some $j\in [n]$ such that $v_j\in e_i$ and $z_j\geq 1$ if and only if the set $v(\bz)=\{v_j \mid z_j\geq 1\}$ is a vertex cover for $H$. The linear program \eqref{eq:tauLP} seeks to minimize the cardinality of the vertex cover $v(\bz)$, in accordance to the definition of the transversal number.
\end{proof}

An {\em edge packing} or {\em matching} of a hypergraph $H=(V,E)$ is a subset of disjoint edges, i.e. $D\subseteq E$ such that no two elements of $D$ share a vertex. Since the edges of $H$ are in bijection with  the minimal monomial generators of  the edge ideal $I(H)$ and two edges are disjoint if and only if the monomials representing them in $I(H)$ are coprime, we have  the following description for edge packing in algebraic terms:
\begin{rem}
\label{rem:edgepacking}
Edge packings of a hypergraph $H$ are in bijection with subsets of the minimal monomial generators of  $I(H)$ in which the elements are pairwise coprime. 
\end{rem}

 Maximal and maximum edge packings are defined to be the edge packings that are maximal with respect to containment and to cardinality, respectively. The size of a maximum edge packing is called the {\em packing number} of $H$, denoted $\pi(H)$.
The packing number of a hypergraph is also the solution to an integer optimization problem with constraints given by the incidence matrix \eqref{eq:matB}, which we now describe. 
\begin{lem}
\label{packingLP}
The packing number of a hypergraph $\pi(H)$ is the optimum value of the following integer linear program
\begin{equation}
\label{eq:piLP}
\begin{tabular}{rl}
\emph{maximize} & $y_1+\cdots+y_t $\\
\emph{subject to} & $B{\bf y} \leq {\bf 1}$ \emph{and}  ${\bf y}\in \N^t$.
\end{tabular}
\end{equation}
\end{lem}
\begin{proof}
A vector $\by\in \N^t$ satisfies the inequality $B{\bf y} \leq {\bf 1}$ if and only if for each vertex $v_i\in V$ there at most one $j\in [n]$ such that $v_i\in e_j$ and $y_j\geq 1$ if and only if the set $e(\by)=\{e_i \mid y_i\geq 1\}$ is a packing for $H$. The linear program \eqref{eq:piLP} seeks to maximize the cardinality of the packing $e(\by)$, in accordance to the definition of the packing number.
\end{proof}

Solving integer optimization problems is much harder than solving linear optimization problems in $\R^n$ because the simplex algorithm solves the latter problem efficiently, while there are no efficient algorithms to solve the former. Therefore a standard practice is to consider the real relaxation of an integer program. The relaxations of the integer programs in  \eqref{eq:tauLP} and \eqref{eq:piLP} are described below. This follows a well established trend to study fractional invariants of combinatorial structures; see \cite{Fractional} for an overview of this method.

\begin{defn}
\label{fractransversal}
Define the {\em fractional transversal number} $\tau_f(H)$, of a hypergraph $H$ to be the optimum solution for the following linear program 
\begin{equation}
\label{eq:taufLP}
\begin{tabular}{rl}
\emph{minimize} & $z_1+\cdots+z_n $\\
\emph{subject to} & $B^T{\bf z} \geq {\bf 1}$ \emph{and}  ${\bf z} \geq {\bf 0}$.
\end{tabular}
\end{equation}
\end{defn}

The feasible set of the linear program above is termed the {\em set covering polyhedron} in \cite{MonAlgBook}. This polyhedron is defined as follows 
\begin{equation}
\label{eq:Q}
\mathcal{Q}(H)=\{\bz\in \R^n \mid B^T\bz \geq {\bf 1},  {\bf z} \geq {\bf 0}\}.
\end{equation}
If $\mathcal{Q}(H)$ is an integer polyhedron, meaning that its vertices have integer coordinates, then the hypergraph $H$ is called {\em Fulkersonian} or {\em ideal}. Any Fulkersonian hypergraph $H$ satisfies $\tau_f(H)=\tau(H)$. This class of hypergraphs is analyzed from an algebraic perspective in \cite{TrungFulkersonian, VillarrealFulkersonian} . 

\begin{defn}
\label{fracpacking}
Define the {\em fractional packing number} of a hypergraph $H$, $\pi_f(H)$, to be the optimum solution for the relaxation of packing problem, namely
\begin{equation}
\label{eq:pifLP}
\begin{tabular}{rl}
\emph{maximize} & $y_1+\cdots+y_t $\\
\emph{subject to} & $B{\bf y} \leq {\bf 1}$ \emph{and}  ${\bf y} \geq {\bf 0}$.
\end{tabular}
\end{equation}
\end{defn}

An important tool in linear programming is linear program duality. This is exemplified by the linear programs \eqref{eq:taufLP} and \eqref{eq:pifLP}, which are dual to each other.  A core aspect of linear optimization is that (real) dual linear programs have the same optimum value, hence  there is an equality $\pi_f(H)=\tau_f(H)$. In fact, based on linear programing duality, \Cref{lem:transversalopt}, and \Cref{packingLP} we deduce the following inequalities 
\begin{equation}
\label{eq:duality}
\pi(H)\leq \pi_f(H)=\tau_f(H) \leq \tau(H).
\end{equation}

It is natural to ask under what circumstances there is equality among the four invariants involved in equation \eqref{eq:duality}.
In combinatorial optimization one considers more generally pairs of dual integer programs with arbitrary objective function. When these pairs have equal optimum values the hypergraph is said to satisfy the {\em max-flow min-cut property}.

We note that the equality $\pi(H)=\tau(H)$ is equivalent to asking for the edge ideal $I(H)$ to be K\"onig cf. \Cref{def:packed}.

\begin{rem}
\label{rem:Konighypergraph}
A hypergraph $H$ satisfies the equality $\pi(H)=\tau(H)$ if and only if $I(H)$ is K\"onig. Indeed, by \Cref{rem:edgepacking}, $\pi(H)$ is the cardinality of the largest set of pairwise coprime monomials among the generators of $I$, whereas by \Cref{vertexcovers} $\tau(H)=\het(I(H))$.
\end{rem}

The following celebrated theorem of K\"onig and Egerv\'ary provides a context in which the equality $\pi(H)=\tau(H)$ is achieved for graphs. Together with the preceding remark, this shows that edge ideals of bipartite hypergraphs are K\"onig.

\begin{thm}[K\"onig, Egerv\'ary -- see e.g.~{\cite[Theorem 2.1.1, p.30]{Diestel}}]
\label{thm:Konig}
In any bipartite graph $G$, the number of edges in a maximum matching equals the number of vertices in a minimum vertex cover, i.e $\tau(G)=\pi(G)$.
\end{thm}

\subsection{The packing problem as a combinatorial optimization problem}
In this section we reformulate the packing problem \Cref{conj:packing} in terms of the linear optimization invariants of hypergraphs introduced above.

\Cref{rem:Konighypergraph} suggests the following definition:
\begin{defn} \label{defn:konig}
A hypergraph $H$ is {\em K\"onig} if it satisfies the equality $\tau(H)=\pi(H)$.
\end{defn}

Following the convention in \cite{Corn, francisco2013powers}, we define two operations on hypergraphs to get smaller hypergraphs. 

\begin{defn}
\label{defn:deletioncontraction}
A \emph{deletion} in a hypergraph is the removal of a vertex $v$ from the vertex set and the removal of any edges that contain it from the edge set. A \emph{contraction} in a hypergraph is the removal of a vertex from the vertex set and from any edges that contain it.

A {\em minor} of a hypergraph $H=(V,E)$ is a hypergraph obtained through a sequence of deletions and contractions. More precisely, it is a hypergraph 
\[h=(V\setminus (V'\cup V"), \{e \setminus V'' \mid e \in E,\,  e\cap V'= \emptyset\}\]
obtained by fixing disjoint (possibly empty) sets $V', V''\subseteq V$, deleting all vertices in $V'$ and contracting all vertices in $V''$.
\end{defn}

We translate \Cref{defn:deletioncontraction} into algebraic language as follows.

\begin{lem} \label{lem:minors}
If $h$ is a minor of $H$ then $I(h)$ is obtained from $I(H)$ by setting the variables corresponding to $v\in V'$ equal to 0 and the variables corresponding to $v''\in V''$ equal to 1.
\end{lem}
\begin{proof}
This follows from the description of the edge set of the minor $h$ in \Cref{defn:deletioncontraction}.
\end{proof}

Based on \Cref{defn:konig} and \Cref{lem:minors}, we can translate the packing property of square-free monomial ideals \Cref{def:packed} into an equivalent definition for  hypergraphs.

\begin{defn} \label{defn:packing_property}
A hypergraph $H$ is said to have the {\em packing property} if  every minor $h$ of $H$ is K\"onig, that is, satisfies $\tau(h)=\pi(h)$.
\end{defn}

Finally, we can restate the packing problem in combinatorial language:

\begin{conj}[The packing problem - hypergraph version] \label{conj:packing}
A hypergraph $H$ satisfies $I(H)^{(m)}=I(H)^m$ for all $m\in \N$ if and only if $H$ has the packing property. 
\end{conj}

\section{Linear optimization invariants of monomial ideals}
\label{s:alginv}

In this section we introduce some algebraic invariants of monomial ideals which can be realized as solutions of linear optimization problems and we expand upon their relationship to the combinatorial optimization invariants from the previous section.

\subsection{Convex bodies associated to monomial ideals}

For a homogeneous ideal $I$ the initial degree, denoted $\alpha(I)$, is the least degree of a non zero element of $I$.
We show that the initial degree of a monomial ideal can be expressed as the solution of a linear program.
For this, we first define the feasible region of the program.

 \begin{defn}
 \label{def:NP}
 The {\em Newton polyhedron} of a monomial ideal $I$ is the convex hull of the exponent vectors of all monomials in $I$, namely
 \[
 NP(I)=\conv \{(a_1,\ldots, a_n)\in \N^n \mid x_1^{a_1}\cdots x_n^{a_n}\in I\}.
 \]
 \end{defn}
  Newton polyhedra of monomial ideals $I$ are integer (or lattice) polyhedra, meaning that their vertices have integer coordinates. Indeed, the vertices of $NP(I)$ are the exponent vectors for a subset of the minimal generators of $I$. 
  
With this notation, the initial degree of a monomial ideal $I$ can be expressed as the solution of a linear program as follows.

\begin{lem}
\label{lem:alphaNP}
If $I$ is a monomial ideal then the initial degree $\alpha(I)$ is the solution of the following linear program
\begin{equation}
\label{eq:alphaLP}
\begin{tabular}{rl}
\emph{minimize} & $a_1+\cdots+a_n $\\
\emph{subject to} & $\ba=(a_1,\ldots, a_n)\in NP(I)$.
\end{tabular}
\end{equation}
\end{lem}
\begin{proof}
This follows because the optimal solution is attained at a vertex of $NP(I)$ and, as remarked above, the vertices of $NP(I)$ correspond to a subset of the minimal generators of $I$. Thus the optimum value of the linear program \eqref{eq:alphaLP} corresponds to a minimal generator of $I$ of least degree.
\end{proof}

While the vertices of the Newton polyhedron are easy to understand, the dual description in terms of bounding inequalities is often difficult to come by. Below we describe a different polyhedron obtained from a square-free monomial ideal which has the advantage that its bounding inequalities can be read off the prime decomposition of the ideal.

 \begin{defn}
 \label{def:SP}
 The {\em symbolic polyhedron} of a square-free monomial ideal $I$ with prime decomposition $I=P_1\cap \cdots \cap P_s$ such that $P_j=(x_{j_1}, \cdots, x_{j_{s_j}})$ for $j\in[s]$  is defined to be the intersection of the Newton polyhedra of the prime components
 \[
 SP(I)=NP(P_1)\cap \cdots \cap NP(P_s)
 \]
 Equivalently, $\by=(y_1,\ldots, y_n)\in \R^n$ is a point in $SP(I)$ if and only if it satisfies 
 \[
 \begin{cases}
 y_{j_1}+\cdots +y_{j_{s_j}}\geq 1 & \text{for } 1\leq j\leq s \\
 y_i\geq 0& \text{for } 1\leq i\leq n.
 \end{cases}
 \]
 \end{defn}
In contrast to the Newton polyhedron, the symbolic polyhedron is a rational polyhedron, meaning that its vertices have rational coordinates. For applications of the symbolic polyhedron, including relationships to combinatorics, see \cite{bocci2016waldschmidt, cooper2017symbolic}.

The following result elucidates the relationship between the Newton and symbolic polyhedra of a monomial ideal.

\begin{prop}
\label{prop:square-free_lattice_points}
Let $I$ be a square-free monomial ideal. Then, there is a containment $NP(I)\subseteq SP(I)$ and the two polyhedra have the same lattice points, that is, 
\[
NP(I)\cap \N^n =SP(I)\cap \N^n.
\]
\end{prop}

\begin{proof}
Let $I$ be a square-free monomial ideal, with decomposition into prime monomial ideals given by $I=P_1\cap\cdots\cap P_s$ for some $s\in\N$. The containment $NP(I)\subseteq SP(I)$ follows from the considering the containments $I\subseteq P_i$ which yield $NP(I)\subseteq NP(P_i)$ for $i\in [s]$. Therefore we conclude $NP(I)\subseteq \bigcap_{i=1}^s NP(P_i)=SP(I)$.

The previous containment implies $NP(I)\cap \N^n \subseteq SP(I)\cap \N^n$. Let $\ba=(a_1,\ldots,a_n)\in SP(I)\cap \N^n$ be a lattice point in $SP(I)$. It follows that for all $i\in [s]$ we have $\ba\in NP(P_i)$. It is well known that the lattice points in the Newton polyhedron of a monomial ideal correspond to monomials in the integral closure of the ideal  \cite[Proposition 1.6]{SH}, hence 
$\bx^{\ba}:=x_1^{a_1}\cdots x_n^{a_n}\in \ov{P_i}$, where $\ov{P_i}$ denotes the integral closure of $P_i$. Since monomial prime ideals are integrally closed, we conclude that $\bx^{\ba}\in P_i$ for all $i\in[s]$, thus $\bx^{\ba}\in I$ and $\ba\in NP(I)\cap \N^n$, as desired.
\end{proof}

We can now give an alternate description for the initial degree of a square-free monomial ideal. To do this we need to associate a matrix to the prime decomposition of a monomial ideal.

\begin{defn}
\label{def:matA}
For a square-free monomial ideal $I$ with prime decomposition $I=P_1\cap \cdots \cap P_s$ such that $P_j=(x_{j_1}, \cdots, x_{j_{s_j}})$ for $j\in[s]$ we define a $s\times n$ {\em prime decomposition matrix} with entries
\begin{equation}
\label{eq:matA}
A_{i,j} = \begin{cases}
1 & \mbox{if $x_j \in P_i$} \\
0 & \mbox{if $x_j \not\in P_i$.}
\end{cases}
\end{equation}
\end{defn}

In the following statement $\ba=\begin{bmatrix} a_1& \cdots & a_n\end{bmatrix}^T$ denotes a vector in $\R^n$.

\begin{lem}
\label{lem:alphaLP}
If $I$ is a square-free monomial ideal then the initial degree $\alpha(I)$ is the optimal solution of the following equivalent linear programs
\begin{equation}
\label{eq:alphaLP}
\begin{tabular}{rlcrl}
\emph{minimize} & $a_1+\cdots+a_n $ & & \emph{minimize} & $a_1+\cdots+a_n $\\
\emph{subject to} & $\ba \in SP(I)\cap \N^n$ & \qquad \qquad &\emph{subject to} & $A\ba\geq {\bf 1}$ \emph{and} $\ba\in \N^n.$
\end{tabular}
\end{equation}
\end{lem}
\begin{proof}
That $\alpha(I)$ is the optimal solution of the leftmost linear program follows from \Cref{lem:alphaNP}, \Cref{prop:square-free_lattice_points} and the fact that $NP(I)$ is an integer polyhedron. It remains to show the equivalence of the two linear programs in the statement. This reduces to showing they have the same feasible region. Comparing the bounding inequalities for $SP(I)$ provided in \Cref{def:SP} to the inequalities in the rightmost linear program defined using the prime decomposition matrix \eqref{eq:matA} one concludes that they coincide.
\end{proof}

We now turn to the relaxation of the linear program in \eqref{eq:alphaLP}, which yields a fractional version of the initial degree.
It turns out that this algebraic invariant has first appeared in the literature under a different guise, which we now recall.

\begin{defn}[\cite{BoH}]
\label{def:Waldschmidt}
The {\em Waldschmidt constant} of a homogeneous ideal $I$ is the value of the following limit
\[
\widehat{\alpha}(I)=\lim_{m\to\infty} \frac{\alpha(I^{(m)})}{m}.
\]
\end{defn}
It turns out that the sequence $\{\alpha(I^{(m)})\}_{m\in \N}$ is subadditive as shown by  the containments $I^{(m)}I^{(m')}\subseteq I^{(m+m')}$ for all $m,m'\in\N$.  Farkas's lemma \cite{Farkas} thus applies to show that  the limit in \Cref{def:Waldschmidt} exists and is equal to the infimum of the respective sequence.

The following theorem shows that the Waldschmidt constant of a square-free monomial ideal is the optimum value of the relaxation of the linear program \eqref{eq:alphaLP}.

\begin{thm}[{\cite[Corollary 6.3]{cooper2017symbolic}}]
      \label{thm:WaldschmidtLP}
For a square-free monomial ideal $I$ with prime decomposition matrix $A$ as in \eqref{eq:matA},  the Waldschmidt constant $\widehat{\alpha}(I)$ is the optimum value of the following equivalent linear programs
\begin{equation}
\label{eq:WaldschmidtLP}
\begin{tabular}{rlcrl}
\emph{minimize} & $a_1+\cdots+a_n $ & & \emph{minimize} & $a_1+\cdots+a_n $\\
\emph{subject to} & $\ba \in SP(I)$ & \qquad \qquad &\emph{subject to} & $A\ba\geq {\bf 1}$ \emph{and} $\ba\geq {\bf 0}.$
\end{tabular}
\end{equation}
\end{thm}

\subsection{Alexander duality}

To relate the algebraic invariants for monomial ideals introduced in \cref{s:alginv} to the combinatorial invariants for hypergraphs encountered in \cref{s:hypinv} it is convenient to introduce the notion of Alexander duality.

\begin{defn}
\label{dual_graph_ideal}
    Let $I$ be a square-free monomial ideal with prime decomposition $I = P_1 \cap \cdots \cap P_s$, where $P_j = ( x_{j_1}, \cdots, x_{j_{s_j}})$ for $j\in[s]$. The \emph{Alexander dual} of $I$ is  the square-free monomial ideal
    \[ I^{\lor} = ( x_{j_1}x_{j_2} \cdots x_{j_{s_j}} \mid j \in [s]). \] 

If $H$ is the hypergraph with edge ideal $I=I(H)$ we define the \emph{dual hypergraph} of $H$ as the hypergraph $H^\vee$ whose edge ideal is $I(H)^\vee$, i.e., $I(H^\vee) = I(H)^\vee$.  In the combinatorial optimization literature $H^\vee$ is called the {\em blocker} of $H$.
\end{defn}

The importance of Alexander duality in our setting is that it interchanges the prime
decomposition matrix \eqref{eq:matA} and the (transpose of the) incidence matrix \eqref{eq:matB}. The following observation arises from comparing \Cref{def:matA}, \Cref{def:matB}, and \Cref{dual_graph_ideal}.

\begin{lem}
If $I$ is a square-free monomial ideal, $H$ is the hypergraph satisfying $I=I(H)$, $A$ and $A^\vee$ denote the prime decomposition matrices of $I$ and $I^\vee$ respectively, and $B$ and $B^\vee$ denote the incidence matrices of $H$ and $H^\vee$ respectively, then 
\[
A^\vee =B^T \qquad \text{ and } \qquad B^\vee=A^T.
\]
\end{lem}

This simple observation shows how the algebraic invariants of monomial ideals relate to combinatorial invariants of hypergraphs.

\begin{cor}
If $I$ is a square-free monomial ideal and $H$ is the hypergraph satisfying $I=I(H)$, then the initial degree and Waldschmidt constant of $I$ can be expressed in terms of the (fractional) transversal number of the dual hypergraph $H^\vee$ as follows:
\[
\alpha(I)=\tau(H^\vee), \quad \widehat{\alpha}(I)=\tau_f(H^\vee) \qquad \text{ and } \qquad \tau(H)=\alpha(I^\vee),\quad \tau_f(H)=\widehat{\alpha}(I^\vee).
\]
\end{cor}

We now turn to convex geometric relationships between the symbolic polyhedron of a square-free monomial ideal, and the set covering polyhedra $\mathcal{Q}(H)$ and $\mathcal{Q}(H^\vee)$ for the corresponding hypergraph and its dual.

\begin{cor}
\label{cor:SPintegral}
Let $I$ be a square-free monomial ideal and let $H$ be the hypergraph satisfying $I=I(H)$. The following are equivalent 
\begin{enumerate}
\item the symbolic polyhedron of $I$ is an integer polyhedron, 
\item the hypergraph $H$ is Fulkersonian, 
\item the dual hypergraph $H^\vee$ is Fulkersonian,
\item  the symbolic polyhedron of $I^\vee$ is an integer polyhedron.
\end{enumerate}
\end{cor}
\begin{proof}
From \Cref{def:SP}, equation \eqref{eq:Q}, and the previous lemma, it follows that $SP(I)=\mathcal{Q}(H^\vee)$. The latter is an integer polyhedron if and only if $H^\vee$ is Fulkersonian, establishing the equivalence of (1) and (2). The equivalence of (2) and (3) follows from \cite[Corollary, p. 210]{Berge}. The equivalence of (3) and (4) follows from the equivalence of (1) and (2) by duality.
\end{proof}

Using \Cref{lem:minors} one can express contraction and deletion as dual operations through the lens of Alexander duality.
\begin{lem}
\label{lem:ADminors}
If $I$ is a square-free monomial ideal, $H$ is the hypergraph satisfying $I=I(H)$, $V', V''$ are subsets of the vertex set of $H$, and $H'$ and $H''$ are the minors of $H$ obtained by deleting $V'$ and contracting $V''$ respectively then 
\begin{enumerate}
\item $(I(H'))^\vee=(I|_{\{x_v'=0\mid v'\in V'\}})^\vee=I^\vee|_{\{x_v'=1\mid v'\in V'\}}$ and
\item $(I(H''))^\vee=I^\vee|_{\{x_v''=0\mid v''\in V''\}}=(I|_{\{x_v''=1\mid v''\in V''\}})^\vee$.
\end{enumerate}
In particular, $(H')^\vee$ is obtained from $H^\vee$ by contracting $V'$ and $(H'')^\vee$ is obtained from $H^\vee$ by deleting $V''$.
\end{lem}
\begin{proof}
\Cref{lem:minors} yields $I(H')=I|_{\{x_v'=0\mid v'\in V'\}}$ and $I(H'')=I|_{\{x_v''=1\mid v''\in V''\}}$. Suppose $I=P_1\cap \cdots\cap P_s$ is the irredundant prime decomposition of $I$ and that $P_i\cap V''=\emptyset$ if and only if $i\in[t]$. Set $P'_i=(x_v \mid x_v\in P_i\setminus V')$, $m_i=\prod_{x_v\in P_i}x_v$ and $m'_i=\prod_{x_v\in P'_i}x_v=m_i|_{\{x_v'=1\mid v'\in V'\}}$ . Then a prime decomposition of $I(H')$ is
\[
I(H')=\bigcap_{i=1}^s P_i|_{\{x_v'=0\mid v'\in V'\}}=\bigcap_{i=1}^s P'_i,
\]
which yields $I(H')^\vee=(m'_1,\ldots, m'_s)=(m_1,\ldots, m_s)|_{\{x_v'=1\mid v'\in V'\}}=I^\vee|_{\{x_v'=1\mid v'\in V'\}}$.
Moreover, a prime decomposition of $I(H'')$ is
$
I(H'')=\bigcap_{i=1}^t P_i
$
which yields 
\[I(H'')^\vee=(m_1,\ldots, m_t)=(m_1,\ldots, m_s)|_{\{x_v''=0\mid v''\in V''\}}=I^\vee|_{\{x_v''=0\mid v''\in V''\}}.\]

\end{proof}

\section{Consequences of the packing problem}
\label{s:consequences}

We now turn our attention to establishing some consequences of the packing problem. These are recorded in  \Cref{thm:packing_and_SP}, \Cref{conj:packing_and_LP}, and \Cref{conj:packing_and_waldschmidt}, which constitute the main results of this section. Their interconnections are summarized in the following sequence of implications elaborated upon in below
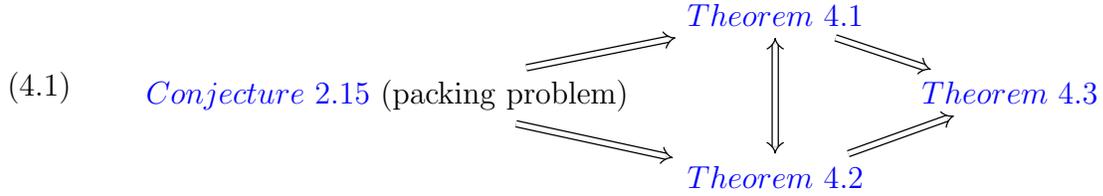
\begin{equation}
\label{eq:implications}
\begin{tikzcd}[row sep=1em,column sep=1em]
& \Cref{thm:packing_and_SP}  \arrow[dd, Leftrightarrow]  \arrow[dr, Rightarrow]  \\
\Cref{conj:packing}\text{ (packing problem)} \arrow[ur, Rightarrow]  \arrow[dr, Rightarrow] &&\Cref{conj:packing_and_waldschmidt}\\
& \Cref{conj:packing_and_LP} \arrow[ur, Rightarrow] 
\end{tikzcd}
\end{equation}
We note that the horizontal implications are non-reversible in \cref{s:reverse}.


We now state three consequences of the packing problem. We will prove their validity in the remainder of the section.

The first is a convex geometric shadow of \Cref{conj:packingalg}.

\begin{thm}
\label{thm:packing_and_SP}
If $I$ is a square-free monomial ideal which satisfies the packing property then $SP(I)=NP(I)$.
\end{thm}

The second is a linear optimization shadow of \Cref{conj:packingalg}.
\begin{thm}
\label{conj:packing_and_LP}
If $I$ is a  square-free monomial ideal which satisfies the packing property and $f(\ba)=c_1a_1+\cdots+c_da_d$ is any linear function with $c_i\geq 0$ for each $i$ then the following two linear programs have equal optimum values:
\begin{equation}
\label{eq:LP}
\begin{tabular}{rlcrl}
\emph{minimize} & $c_1a_1+\cdots+c_da_d$ & \qquad  \qquad & \emph{minimize} & $c_1a_1+\cdots+c_da_d$\\
\emph{subject to} & $\ba\in SP(I)$ & & \emph{subject to} & $\ba\in NP(I)$.
\end{tabular}
\end{equation}
\end{thm}

The restriction $c_i\geq 0$ for each $i$ insures that the optimum values of the linear programs in \eqref{eq:LP} are real numbers. If this is not satisfied, then the optimum values of both programs are $-\infty$. This is  because both the Newton and the symbolic polyhedron are closed under increasing coordinates.

The third consequence is a numerical shadow of \Cref{conj:packingalg}.

\begin{thm} \label{conj:packing_and_waldschmidt}
If $I$ is a  square-free monomial ideal which satisfies the packing property then there is an equality \[\widehat{\alpha}(I)=\alpha(I).\]
\end{thm}

We start by showing that the validity of each of the above theorems follows from the validity of  the packing problem \Cref{conj:packingalg}. For this purpose we recall the following result regarding points with rational coordinates in symbolic polyhedra.

\begin{lem}[{\cite[Proposition 6.1]{cooper2017symbolic}}]
\label{lem:CEHHprop}
Let $I$ be a monomial ideal with symbolic polyhedron $SP(I)$. For any $\ba\in SP(I)\cap\Q^d$ there exists a positive integer $b$ such that $\bx^{m\ba}\in I^{(m)}$ whenever $m$ is divisible by $b$.  
\end{lem}



Armed with this result, we are now ready to prove a general result regarding the equality of the symbolic and Newton polyhedra.

\begin{prop} \label{prop:SP=NP}
Let $I$ be a monomial ideal such that $I^{(n)}=I^n$ for all $n\geq 1$. Then there is an equality $NP(I)=SP(I)$.
\end{prop}
\begin{proof}
To prove the claim, it suffices to show that every vertex of $SP(I)$ lies in $NP(I)$. By convexity of $NP(I)$ and $SP(I)$, this would imply $SP(I)\subseteq NP(I)$. Let $\ba=(a_1,\ldots,a_d)\in SP(I)$ be a vertex. Since the bounding hyperplanes of $SP(I)$ are given by equations with integer coefficients by \Cref{def:SP}, we have that $\ba\in\Q^n$. By \Cref{lem:CEHHprop}, we can choose some positive integer $b$ such that $\bx^{b\ba}\in I^{(b)}= I^b$. This implies that $b\ba\in NP(I^b)=b\cdot NP(I)$, which allows us to conclude that $\ba\in NP(I)$ as required.

The argument above shows that $SP(I)\subseteq NP(I)$. Since the reverse inclusion holds in general (see \Cref{prop:square-free_lattice_points}), we obtain the desired equality.
\end{proof}

The preceding result allows us to show \Cref{conj:packingalg} implies \Cref{thm:packing_and_SP}.

\begin{prop}
\label{prop:packing_implies_polyhedra}
Assume that the assertion of the packing problem, \Cref{conj:packingalg}, is true. Then any square-free monomial ideal $I$ which has the packing property satisfies the equality $NP(I)=SP(I)$.
\end{prop}
\begin{proof}
Let $I$ be a  square-free monomial ideal which satisfies the packing property. Assuming that the statement of the packing problem is true, the hypothesis implies that the equalities $I^{(n)}=I^n$ hold for all $n\in \N$. By \Cref{prop:SP=NP}, it follows that $SP(I)=NP(I)$. 
\end{proof}

We next show the equivalence of \Cref{thm:packing_and_SP} and \Cref{conj:packing_and_LP} connecting  the equality of the Newton and symbolic polyhedra and the equivalence of linear programs with nonnegative coefficients for the objective function.

%
%

\begin{prop}
Let $I$ be a monomial ideal. Then there is an equality of polyhedra $SP(I)=NP(I)$ if and only if the following linear programs have the same optimal solution for all objective functions $f(\ba)=c_1a_1+\cdots+c_na_n$  with $c_i\geq 0$ for each $i\in [n]$ 
\begin{equation}
\label{eq:4.3}
\begin{tabular}{rlcrl}
\emph{minimize} & $c_1a_1+\cdots+c_na_n$ & \qquad  \qquad & \emph{minimize} & $c_1a_1+\cdots+c_na_n$\\
\emph{subject to} & $\ba\in SP(I)$ & & \emph{subject to} & $\ba\in NP(I)$.
\end{tabular}
\end{equation}
\end{prop}
\begin{proof}
The forward implication is clear. We now focus on the converse.

Consider the symbolic polyhedron $SP(I)$ and the Newton polyhedron $NP(I)$ of $I$ and recall from \Cref{prop:square-free_lattice_points} that $NP(I)\subseteq SP(I)$. To show that $NP(I)=SP(I)$, it suffices to show that the vertices of $SP(I)$ are contained in $NP(I)$. 
Let $\bp=(p_1,\ldots,p_d)$ be a vertex of $SP(I)$. Note that $\bp$ is in fact the intersection point of $n$ distinct bounding hyperplanes $H_1,\ldots, H_n$ for $NP(I)$. By \Cref{def:SP}, these bounding hyperplanes have non-negative (with entries  $0$ and $1$) normal vectors $\bc_1,\ldots, \bc_n$ respectively. Taking their arithmetic mean $$\bc=\frac{1}{n}(\bc_1+\ldots+ \bc_n)=(c_1,\ldots,c_n)\in\R_{\geq 0}^n,$$ 
we obtain a hyperplane $H_\bp$ with equation $c_1(a_1-p_1)+\cdots +c_n(a_n-p_n)=0$ which intersects $SP(I)$ only at $\bp$.
This is because the equation of $H_\bp$ is the arithmetic mean of the equations of $H_1,\ldots, H_n$ and for each point of $SP(I)$ other than $\bp$ the result of substituting it into the equations of $H_1,\ldots, H_n$ is always non-negative, with at least one positive value.
Since $H_\bp \cap SP(I)=\{\bp\}$,  the optimal value of the linear program 
\begin{equation*}
\label{eq:mean}
\begin{tabular}{rl}
minimize & $f(\ba)=c_1a_1+\cdots+c_da_d$ \\
subject to & $\ba\in SP(I)$ 
\end{tabular}
\end{equation*}
is attained at $\mathbf{p}$. Since $\bc\geq {\bf 0}$, the hypothesis implies that $f(\bp)$ is the optimal value of the  linear program with objective function $f$ and feasible set  $NP(I)$. Since $f(\ba)>f(\bp)$ for all points $\ba\in SP(I)\setminus\{\bp\}$ and since $NP(I)\subseteq SP(I)$, it follows that $f(\ba)>f(\bp)$ for all points $\ba\in NP(I)\setminus\{\bp\}$.  We conclude that the point $\bp$ must belong to $NP(I)$ in order for $f(\bp)$ to be the minimum value of the second linear program in \eqref{eq:4.3}. Since $\bp$ was an arbitrary vertex of $SP(I)$, it follows that $SP(I)\subseteq NP(I)$, as desired.
\end{proof}

The final implication needed to complete diagram \eqref{eq:implications} is the following.

\begin{lem} \Cref{conj:packing_and_LP} implies \Cref{conj:packing_and_waldschmidt}.
\end{lem}
\begin{proof}
Setting $c_i=1$ for $i\in [n]$ in the linear programs \eqref{eq:LP} yields the linear programs \eqref{eq:WaldschmidtLP} and \eqref{eq:alphaLP} respectively. By \Cref{lem:alphaLP} and \Cref{thm:WaldschmidtLP} the optimal values of these programs are $\alpha(I)$ and $\widehat{\alpha}(I)$. Thus the conclusion of \Cref{conj:packing_and_LP} implies that $\alpha(I)=\widehat{\alpha}(I)$, i.e., the conclusion of \Cref{conj:packing_and_waldschmidt}.
\end{proof}

Finally, we prove \Cref{conj:packing_and_waldschmidt}, \Cref{thm:packing_and_SP},  and \Cref{conj:packing_and_LP} independent of the (not yet established) validity of the packing problem. Due to the implications in diagram \eqref{eq:implications}, it suffices to prove the validity of  \Cref{thm:packing_and_SP}, since this result implies the other two.
Towards this end we use a celebrated result of Lehman.

\begin{thm}[\cite{Lehman}, {\cite[Theorem 1.8]{conforti1990decomposition}}]
\label{thm:Lehman}
If a hypergraph $H$  has the packing property, then the polyhedron $\mathcal{Q}(H)$ defined in equation \eqref{eq:Q} 	is an integer polyhedron.
\end{thm}

We are now ready to prove our main results.

\begin{proof}[Proof of \Cref{thm:packing_and_SP}]

Let $H$ and $H^\vee$ be the hypergraphs determined by $I(H)=I$ and $I(H^\vee)=I^\vee$ respectively. Since $I$ and hence $H$ satisfy the packing property by hypothesis, \Cref{thm:Lehman} yields that the set covering polyhedron $\mathcal{Q}(H)$ is an integer polyhedron, thus $H$ is Fulkersonian. \Cref{cor:SPintegral} now yields that $I$ has an integer symbolic polyhedron. Since the vertices of $SP(I)$ are lattice points, they belong to $NP(I)$ by \Cref{prop:square-free_lattice_points},  thus inducing a containment $SP(I)\subseteq NP(I)$. Since the opposite containment always holds (see \Cref{prop:square-free_lattice_points}) we conclude the desired equality $SP(I)=NP(I)$.
\end{proof}

As previously remarked, the lattice points in the Newton polyhedron $NP(I)$ of any monomial ideal $I$ correspond to monomials in the integral closure $\overline{I}$ of the ideal $I$; see \cite[Proposition 1.6]{SH}). Thus the Newton polyhedra $NP(I)$ and $NP(\overline{I})$ coincide. Thus the equality $SP(I)=NP(I)$ can be rewritten as $SP(I)=NP(\overline{I})$ and  thought of as capturing the equality of symbolic powers and integral closures of powers of $I$. We present an alternate proof of \Cref{thm:packing_and_SP} based on this intuition. For this purpose we recall an alternate description of the symbolic polyhedron from \cite{Waldschmidtpaper}.

\begin{lem}[{\cite[Corollary 3.12]{Waldschmidtpaper}}]
Let $I$ be a monomial ideal. Then the symbolic polyhedron of $I$ can be described as 
\[
SP(I)=\bigcup_{m\geq 1}\frac{1}{m}NP(I^{(m)}).
\]
\end{lem}

\begin{proof}[Alternate proof of \Cref{thm:packing_and_SP}]
As in the previous proof of this result, the hypothesis that $I$ sa\-tis\-fies the packing property implies that $H$ is Fulkersonian. By \cite[Theorem 2.3]{TrungFulkersonian} or \cite[Proposition 3.4]{VillarrealFulkersonian} this guarantees  equality of the symbolic powers and integral closures of ordinary powers, namely $I^{(m)}=\overline{I^m}$ for each $m\geq 1$. Passing to the respective convex bodies yields the identities 
\[NP(I^{(m)})=NP(\overline{I^m})=NP(I^m)=mNP(I)\]
and taking the convex limit of the above family of polyhedra yields the desired equality
\[
SP(I)=\bigcup_{m\geq 1}\frac{1}{m}NP(I^{(m)})=\bigcup_{m\geq 1}\frac{1}{m}\cdot m NP(I)=\bigcup_{m\geq 1} NP(I)=NP(I).
\]
\end{proof}

\section{Further questions and conjectures}
\label{s:reverse}

In this section we consider the implications that our results have on the packing problem. This amounts to reversing the implications in diagram \eqref{eq:implications}. While we show below that in general these implications are not reversible, this line of reasoning leads us to some related conjectures that have a bearing on the packing problem.

\subsection{Uniform hypergraphs}

In \Cref{prop:packing_implies_polyhedra} we showed that the packing problem implies the validity of \Cref{thm:packing_and_SP}. We note that \Cref{thm:packing_and_SP},  {\em not} imply the validity of \Cref{conj:packing} (the packing problem) and neither do \Cref{conj:packing_and_LP} or \Cref{conj:packing_and_waldschmidt}, as illustrated by the following \Cref{ex:DD}. This is closely related to the irreversibility of \Cref{thm:Lehman} also demonstrated by this example. We first learned of \Cref{ex:DD} from \cite[Remark 5.4]{DD}.
In combinatorial optimization the corresponding hypergraph has gained some recognition under the name $\mathcal{Q}_6$, see \cite[Example 14.2.9]{MonAlgBook}. It is a forbidden minor of any hypergraph that satisfies the max-flow min-cut property.

\begin{ex}
 \label{ex:DD}
Consider the square-free monomial ideal 
$$I=(abc,aef,cde,bdf)\subseteq K[a,b,c,d,e,f]$$ 
with prime decomposition $$I=(a,d)\cap(b,e)\cap(c,f)\cap(a,b,c)\cap(a,e,f)\cap(b,d,f)\cap(c,d,e),$$ which implies that $\het(I)=2$. It can be verified using a computer algebra system such as {\em Macaulay2} {\cite{M2}} that there is an equality $NP(I)=SP(I)$ and consequently $SP(I)$ is an integer polyhedron and $\alpha(I)=\widehat{\alpha}(I)$.

However, the ideal $I$ does not satisfy the packing property. In particular, we see that the ideal $I$ itself is not K\"onig as any $\het(I)=2$ monomials in $I$ have non-trivial common divisor. This  ideal fails to satisfy the packing property in a minimal way, since it is only the full ideal $I$ that is not K\"onig. All other minors obtained by setting any number of variables equal to 1 or 0 are K\"onig.
\end{ex}

\begin{defn}
   A hypergraph $H$ is {\em uniform} if every edge of $H$ has the same number of vertices, equivalently if the edge ideal $I(H)$ is equigenerated.   An ideal is {\em equidimensional} if all its associated primes have the same height.
\end{defn}

The two notions above are related by Alexander duality: if $I$ is the edge ideal of a hypergraph $H$, then $I$ is equidimensional if and only if $H^\vee$ is uniform and $H$ is uniform if and only if $I^\vee$ is equidimensional.

To establish \Cref{conj:packingalg} for all uniform hypergraphs, it suffices to prove it for those that have Cohen-Macaulay edge ideals; see \cite[Theorem 3.3]{DRV}. The Cohen-Macaulay property forces the edge ideal to be equidimensional. Hence to establish \Cref{conj:packingalg} for all uniform hypergraphs, it suffices to prove it for uniform hypergraphs that have equidimensional edge ideals. See \Cref{prop:equidim} for a result regarding this family of ideals.

We note that the ideal in  \Cref{ex:DD} is not equidimensional. We do not know whether \Cref{thm:packing_and_SP} and \Cref{conj:packingalg} (the packing problem) are equivalent for {\em equidimensional} ideals. This motivates the following question.

\begin{quest}
\label{q:equidimensional}
Is there an equidimensional square-free monomial ideal $I$ so that there is an equality of polyhedra $SP(I)=NP(I)$ (equivalently, $SP(I)$ is an integer polyhedron), but $I$ is not packed or $I^m\neq I^{(m)}$ for some $n$?
\end{quest}

If the answer to \Cref{q:equidimensional} is negative, then this means \Cref{thm:packing_and_SP} and the packing problem are equivalent for equidimensional ideals. Since we have proven \Cref{thm:packing_and_SP}, this would imply the validity of the packing problem for equidimensional ideals.

\subsection{Partite hypergraphs}

Recall that a graph satisfies the packing property if and only if it is bipartite. In this section we investigate notions of partite hypergraphs and their relationship to Waldschmidt constants of edge ideals. We make a conjecture in this regard, which would provide a new bridge from \Cref{conj:packing} to \Cref{conj:packing_and_waldschmidt}.

\begin{notation}
    Let $A$ be a set, and $k$ be a positive integer. We denote by $\binom{A}{k}$ the set of $k$ element subsets of $A$, i.e.
    $\binom{A}{k} = \{ S \subseteq A : |S| = k \}$. 
\end{notation}

\begin{defn}
    Let $H = (V, E)$ be a hypergraph and let $a$ and $b$ be positive integers with $a \geq b$. A function $f : V \to \binom{[a]}{b}$ is called an $(a : b)$--{\em partition} or $(a : b)$--{\em rainbow coloring} of $H$ if for each edge $e$ and for each color $i$, there is a vertex $v \in e$ such that $i \in f(v)$.
A hypergraph $H$ is said to be  \textit{$(a:b)$--partite} or \textit{$(a:b)$--rainbow colorable} if it has an $(a:b)$--rainbow coloring.

We say that $H$ is {\em $a$-partite} or {\em $a$-colorable} if $H$ has an $(a:1)$--coloring.
\end{defn}

We show below how the property of a graph of being $(a:b)$--partite imposes a lower bound on the Waldschmidt constant of its edge ideal. Towards this end, we first introduce a useful property of partite hypergraphs.

\begin{lem}
 \label{prop:vertex_covers_partite}
    Let $H $ be an $(a:b)$--partite hypergraph. Then there exist  disjoint minimal vertex covers $C_1, \ldots, C_a$ such that no vertex  appears in more than $b$ of the vertex covers. 
\end{lem}

\begin{proof}
    Fix an $(a:b)$--rainbow coloring $f$ of $H= (V, E)$, and consider the $a$ color classes of this coloring, i.e., the sets
    $A_i = \{ v \in V \mid i \in f(v) \}$.
    Since $f$ is an $(a:b)$--rainbow coloring of $H$, every color appears at least once in each edge, and thus the color classes are vertex covers. Thus each $A_i$ contains a minimal vertex cover $C_i$. No vertex appears in more than $b$ of the $C_i$, since if it did, that would mean that there is a vertex with more than $b$ colors, which is a contradiction.
\end{proof}

The proof of the following result  uses the description of the Waldschmidt constant as the solution of a linear optimization problem in \Cref{thm:WaldschmidtLP}.

\begin{prop}
 \label{lem:waldschmidt_lower_bound_hypergraphs}
Let $H$ be a hypergraph  with edge ideal $I = I(H)$ and let $a$ and $b$ be positive integers with $a \geq b$. If $H$ is $(a:b)$--partite, then the inequality $\widehat{\alpha}(I) \geq \frac{a}{b}$ holds.
\end{prop}


\begin{proof}
Let 
$I=P_1\cap \cdots \cap P_s$ be the irredundant prime decomposition of $I$. Since $H$ is $(a:b)$--partite, by \Cref{prop:vertex_covers_partite} there is a set of $a$ minimal vertex covers of $H$ such that no vertex appears in more than $b$ of the covers. Minimal vertex covers of $H$ correspond to associated primes of $I$ by \Cref{vertexcovers}. Without loss of generality let $P_1,\ldots,P_a$ be the primes corresponding to these minimal vertex covers. From \Cref{def:SP}, we have
\begin{equation*} \label{eqn:waldschmidt_lower_bound_hypergraphs_1}
SP(I)=NP(P_1) \cap \cdots \cap NP(P_a) \cap \cdots \cap NP(P_s).
\end{equation*}
where for $P_j=(x_{j_1},\ldots, x_{j_{s_j}})$ one has 
\[ NP(P_j) = \left\{ \by \mid y_{j_1}+\cdots +y_{j_{s_j}}\geq 1, y_{j_1}, \cdots, y_{j_{s_j}} \geq 0\right \}. \] 
In particular, any point $\mathbf{c}=(c_1,\ldots, c_n) \in SP(I)$  satisfies the set of $a$ inequalities:
\[ \sum_{k=1}^{s_j} c_{j_k} \geq 1 \qquad\text{ for } j \in [a].\]
Adding up these inequalities one obtains for $\alpha_i=|\{j\in [a] \mid x_i\in P_j\}|$ the equation
\[ \sum_{i=1}^n \alpha_i c_i \geq a.\]
Since no vertex appears in more than $b$ of the $a$ minimal vertex covers corresponding to $P_1, \ldots, P_a$, we have that $\alpha_i \leq b$ for each $i$, leading to the inequality
\[ \sum_{i=1}^n b c_i \geq \sum_{i=1}^n \alpha_i c_i\geq a. \]
It follows that $c_1 + c_2 + \cdots + c_n \geq \frac{a}{b}$. As $\bc$ was arbitrary, minimizing the sum of coordinates over $SP(I)$ as in \Cref{thm:WaldschmidtLP}, we reach the conclusion $\widehat{\alpha}(I) \geq \frac{a}{b}$.
\end{proof}

We have shown in \cref{s:consequences} that the validity of the packing problem for an ideal $I$ implies the equality $\widehat{\alpha}(I)=\alpha(I)$. In  general, for arbitrary ideals $I$, the inequality $\widehat{\alpha}(I)\leq \alpha(I)$ holds thus to have equality it suffices to show the converse inequality. Our next conjecture is inspired by  the lower bound given by \Cref{lem:waldschmidt_lower_bound_hypergraphs}. We conjecture that in the presence of the packing property the inequality in \Cref{lem:waldschmidt_lower_bound_hypergraphs} yields an inequality $\widehat{\alpha}(I)\geq \alpha(I)$ by means of a specific partite structure on the corresponding hypergraph.

\begin{conj}
\label{conj:coloring}
    Let $H$ be a hypergraph with the packing property, and let $I=I(H)$ be its edge ideal. Then $H$ is $(\alpha(I) \het(I) : \het(I))$--partite.
\end{conj}

\subsection{Packing and Alexander duality}
Satisfying the packing property need not be preserved under taking the Alexander dual. In fact the Alexander dual of an ideal $I$ that satisfies the packing property can fail to satisfy this property, even in the case when $I$ is equidimensional and hence  the  hypergraph corresponding to $I^\vee$ is uniform. 
 
\begin{ex} 
\label{ex:the_counterexample}
Consider \Cref{ex:DD} again where we introduced the ideal
\[I = (abc,aef,cde,bdf) = (a,d)\cap(b,e)\cap(c,f)\cap(a,b,c)\cap(a,e,f)\cap(b,d,f)\cap(c,d,e).\]
Its Alexander dual is
 \[I^\vee = (ad,be,cf,abc,aef,bdf,cde) = (a,b,c) \cap (a,e,f) \cap (c,d,e) \cap (b,d,f).\] 
Using the computer algebra system Macaulay2 {\cite{M2}} equipped with the package {\em Sym\-bo\-licPowers} \cite{JSAG}, one can check that $I$ does not satisfy the packing property, while $I^\vee$ does. Notice also that the ideal $I^\vee$ is equidimensional, while $I$ is not equidimensional. 
\end{ex}

This example leads to the following task.

\begin{quest}
\label{quest:dual-packing}
 Give an algebraic or combinatorial description for the class of square-free monomial ideals $I$ such that $I$ satisfies the packing property if and only if $I^\vee$ does. 
    \end{quest}
    
    In \Cref{prop:equidim} we  give a partial answer by showing that equidimensional ideals $I$ which have equidimensional duals $I^\vee$ are part of the class of ideals singled out in \Cref{quest:dual-packing}.
Note that the family of equidimensional ideals $I$ such that $I^\vee$ is also equidimensional corresponds bijectively to uniform hypergraphs $H$ which have a uniform blocker  $H^\vee$.

We begin our inquiry by studying properties of uniform hypergraphs which have the packing property.

\begin{thm}
\label{prop:uniformpacking}
Let $H$ be a uniform hypergraph which satisfies the packing property and let $I=I(H)$. Then
\begin{enumerate}
\item $H$ has an exact cover, meaning a vertex cover which meets every edge in exactly one vertex,
\item $H$ is $\alpha(I)$--partite,
\item $I^\vee$ and $H^\vee$ satisfy the packing property.
\end{enumerate}
\end{thm}
\begin{proof}
Because $H$ satisfies the packing property, the polyhedron $\mathcal{Q}(H)$ is integral by \Cref{thm:Lehman}. A uniform hypergraph with integral set covering polyhedron has an exact cover by \cite[Lemma 5.6 and Proposition 3.13]{GRV}, see also \cite[Lemma 14.4.1]{MonAlgBook}.

We now show that $H$ is $\alpha(I)$--partite by induction on $\alpha(I)$. If $\alpha(I)=1$ the claim follows since every hypergraph is 1-partite. Otherwise, let $C$ be an exact cover of $H$ and denote $h$ the minor of $H$ obtained by deleting the vertices in $C$. The monomial ideal $I(h)$ satisfies the packing property and $\alpha(I(h))=\alpha(I)-1$, therefore $h$ is $(\alpha(I)-1)$--partite by induction. Extending the coloring of $h$ to $H$ by coloring the vertices in $C$ with a new color, yields that $H$ is $\alpha(I)$--partite, as desired.

Let $A_i$ denote the set of vertices colored by color $i\in[\alpha(I)]$. Then each set $A_i$ is a vertex cover and thus contains a  minimal vertex cover $C_i$ of $H$. Since the sets $A_i$ are disjoint, so are the minimal vertex covers $C_i$. Consequently, $I^\vee$ contains a regular sequence of $\alpha(I)$ monomial generators $m_i=\prod_{v\in C_i} x_v$. Since $\het(I^\vee)=\alpha(I)$ by definition, it follows that $I^\vee$ is K\"onig according to \Cref{def:packed}. 

To show that $I^\vee$ has the packing property, consider subsets $V', V''$ of the vertices of $H$. Since $H$ is uniform, the minor $H'$ of $H$ obtained by deleting the vertices in $V'$ is a uniform hypergraph.  Its edge ideal is $I(H')=I|_{\{x_v'=0\mid v'\in V'\}}$. 
The identity from \Cref{lem:ADminors}
\[
(I(H'))^\vee=(I|_{\{x_v'=0\mid v'\in V'\}})^\vee=I^\vee|_{\{x_v'=1\mid v'\in V'\}}\]
reveals that the ideal $I^\vee|_{\{x_v'=1\mid v'\in V'\}}$ is K\"onig, because we have shown above that the duals of uniform 
hypergraphs that satisfy the packing property are K\"onig.

Now consider the minor $H''$ obtained from $H$ by contracting the vertices in $V''$. This need not be a uniform hypergraph.  Its edge ideal is $I(H'')=I|_{\{x_v''=1\mid v''\in V''\}}$. Let \[v=|\{i\mid C_i\subseteq V''\}|.\]
We claim that $H''$ is $(\alpha(I)-v)$--partite. Indeed, $H''$ can be colored with $\alpha(I)-v$ colors using the color assignments $v\mapsto i$ if $v\in C_i\setminus V''$. Since the color sets are disjoint, the monomials $m_i=\prod_{v\in C_i}x_v$ form a regular sequence in $I(H'')^\vee$ of length $\alpha(I)-v$.

Consider the ideal discussed in \Cref{lem:ADminors}
\[
I'':=(I(H''))^\vee=(I|_{\{x_v''=1\mid v''\in V''\}})^\vee=I^\vee|_{\{x_v''=0\mid v''\in V''\}}
\]
and let $\het(I'')=\alpha(I(H''))=u$. 
To establish the claim that $I''$ is K\"onig, it suffices to show that $u\leq \alpha(I)-v$. (In fact this will force $u=\alpha(I)-v$.) Note that the inequality $u\leq \alpha(I)-v$ is equivalent to $v\leq \alpha(I)-\alpha(I(H''))$. Let $m''$ be a minimal generator of $I(H'')$ with $\deg(m'')=\alpha(I(H'')$. Then there is a minimal generator $m$ of $I$ such that $\deg(m)=\alpha(I)$ and $m/m''$ is a product of $\alpha(I)-\alpha(I(H''))$ variables corresponding to vertices in $V''$. Since the sets $C_i$ are vertex covers, each $C_i\subseteq V''$ must contain at least one vertex corresponding to a variable dividing the monomial $m/m''$. Since the  sets $C_i$ are disjoint, this observation yields the desired conclusion $v\leq \alpha(I)-\alpha(I(H''))$.
\end{proof}

\begin{cor}
 \label{prop:equidim}
Let $I$ be an equidimensional square-free monomial ideals $I$ such that $I^\vee$ is also equidimensional.  Then $I$ and  $I^\vee$ satisfy the packing property simultaneously, that is, $I$ satisfies the packing property if and only if $I^\vee$ does.
\end{cor} 
\begin{proof}
The hypothesis implies that the hypergraphs $H$ and $H^\vee$ corresponding to $I$ and $I^\vee$ respectively are uniform. Suppose that any one of $I$ or $I^\vee$ satisfies the packing property. Then part (5) of \Cref{prop:uniformpacking} shows that the dual ideal satisfies the packing property as well.
\end{proof}

{\bf Data availability statement: }
 Data sharing not applicable to this article as no datasets were generated or analysed during the current study.

\vspace{1em}

\bibliographystyle{amsalpha}
\bibliography{references}


\end{document}